\documentclass[12pt, a4paper]{article}
\usepackage{amsmath,amsfonts,amssymb,amsthm,array}
\newcommand{\ds}{\displaystyle}
\newcommand{\NN}{\mathbb{N}}
\newcommand{\RR}{\mathbb{R}}
\renewcommand{\ss}{\mathbf{s}}
\newcommand{\kk}{\mathbf{k}}
\newcommand{\ee}{\mathbf{e}}
\newcommand{\xx}{\mathbf{x}}
\newcommand{\yy}{\mathbf{y}}

\newtheorem{thm}{Theorem}

\begin{document}

\title{A closed formula for the derivatives of $\ds e^{f(x)}$}
\author{Konstantinos Drakakis\footnote{The author holds a Diploma in Electrical and Computer Engineering from NTUA, Athens, Greece, and a Ph.D. in Applied and Computational Mathematics from Princeton University, NJ, USA.}\\School of Mathematics\\University of Edinburgh}
\date{20/6/2005}
\maketitle

\abstract{We give a closed formula for the derivative of arbitrary order of the function $\ds g(x)=\exp(f(x))$.}

\section{Introduction}
The computation of derivatives of functions of the form $\exp(f(x))$ are very useful in analysis, when one is interested in boundedness properties of such functions, and especially in financial mathematics and probability, where this particular form of functions arises naturally. In this paper we will give a closed form for $[\exp(f(n))]^{(n)},\ n\in\NN$. 

\section{Notation}
Let us denote by $\NN^\infty$ the collection of sequences of natural numbers, with the additional restriction that only a finite number of the elements of the sequence can be positive:
\[\NN^\infty=\left\{\kk=\{k_i\}_{i=1}^\infty\bigg|\forall i\in\NN^*, k_i\in\NN \wedge \sum_{i=0}^\infty k_i<\infty\right\}\]
In other words, $\exists n\in\NN: \forall i>n,k_i=0$; in such cases we will also be writing $\kk=(k_1,k_2,\ldots,k_n,0,\ldots)=(k_1,k_2,\ldots,k_n)$, i.e. we will use vector notation. We also define the constant $\ss=(1,2,3,\ldots)$, i.e. $\forall i\in\NN^*, s_i=i$. 

We will write $\ds \kk!=\prod_{i=1}^\infty k_i!=\prod_{i=1}^n k_i!$, and also $\ds \left(f^{(\ss)}\right)^\kk=\prod_{i=1}^\infty \left(f^{(i)}\right)^{k_i}$.

\section{The general formula: a sketch}
We claim that the derivative of order $n$ of the function given is of the form:
\begin{equation}\left(e^{f(x)}\right)^{(n)}=e^{f(x)}\sum_{\{\kk|\kk\in\NN^\infty,\kk\cdot\ss=n\}} a_\kk^n \left(f^{(\ss)}\right)^\kk \label{eq}\end{equation}
We will proceed now to prove that this sum extends over the $\kk$s we claim it does, and also find the exact value of the coefficients $a_\kk^n\in \RR$. Notice that the $n$ in this notation is superfluous, as it can be inferred by $\kk$, but it will prove convenient later, as a marker of the order of the derivative the coefficient belongs to. 

\section{The proof}

\begin{thm}
\emph{(\ref{eq})} holds if $\ds a_\kk^{n+1}=a_{\kk-\ee_1}^n+\sum_{i=1}^n (k_i+1)a_{\kk+\ee_i-\ee_{i+1}}^n$, where the $\ee_i,\ i\in\NN$ denote the sequences in $\NN^\infty$ whose elements are all 0, except the $i$th element which is 1. 
\end{thm}

\begin{proof}
For $n=1$ the formula is true, as $\left(e^{f(x)}\right)'=e^{f(x)}f'(x)$, so that the summation extends indeed over all $\kk$s such that $\kk\cdot\ss=1$, i.e. only on $\kk=(1,0,0,\ldots)=1$, and $a^1_1=1$. 

Assume now that the formula is true for $n$; we need to show it true for $n+1$. But the formula for $n+1$ can be obtained by the formula for $n$ by differentiation:
\begin{multline*}
\left(e^{f(x)}\sum_{\{\kk|\kk\in\NN^\infty,\kk\cdot\ss=n\}} a_\kk^n \left(f^{(\ss)}\right)^\kk\right)'=\\=
e^{f(x)}\sum_{\{\kk|\kk\in\NN^\infty,\kk\cdot\ss=n\}} a_\kk^n f'\left(f^{(\ss)}\right)^\kk+
e^{f(x)}\sum_{\{\kk|\kk\in\NN^\infty,\kk\cdot\ss=n\}} a_\kk^n \left(\left(f^{(\ss)}\right)^\kk\right)'=\\
=e^{f(x)}\sum_{\{\kk|\kk\in\NN^\infty,\kk\cdot\ss=n\}} \left(a_\kk^n \left(f^{(\ss)}\right)^{\kk+\ee_1}+a_\kk^n 
\sum_{i=1}^n k_if^{(i+1)} \left(f^{(\ss)}\right)^{\kk-\ee_i}\right)=\\=
e^{f(x)}\sum_{\{\kk|\kk\in\NN^\infty,\kk\cdot\ss=n\}} \left(a_\kk^n \left(f^{(\ss)}\right)^{\kk+\ee_1}+a_\kk^n 
\sum_{i=1}^n k_i \left(f^{(\ss)}\right)^{\kk-\ee_i+\ee_{i+1}}\right)=\ldots
\end{multline*}
Observe here that $\kk+\ee_1$ and $\kk-\ee_i+\ee_{i+1},\ i=1,\ldots,n$ are solutions of the equation $\ss\cdot\xx=n+1$. On the other hand, if $\xx$ is a solution of this equation, it has at least one positive entry, hence at least one of the vectors $\xx-\ee_1$, $\xx+\ee_i-\ee_{i+1},\ i=1,\ldots,n$ will belong in $\NN^\infty$. But as all such vectors that are in $\NN^\infty$ are solutions of the equation $\ss\cdot\yy=n$, and as, by our inductive assumption, all solutions of this latter equation in $\NN^\infty$ appear in the summation of (\ref{eq}), there will exist a $\kk\in\NN^\infty$ so that either $\xx=\kk+\ee_1$ or $\xx=\kk-\ee_i+\ee_{i+1}$ for some $i=1,\ldots,n$. In other words, the summation in the last formula of the above derivation extends on all $\{\kk|\kk\in\NN^\infty,\kk\cdot\ss=n+1\}$. So, we can go on with the derivation, by rearranging terms:
\begin{multline*}
\ldots
=e^{f(x)}\sum_{\{\kk|\kk\in\NN^\infty,\kk\cdot\ss=n+1\}} \left[a_{\kk-\ee_1}^n+\sum_{i=1}^n a_{\kk+\ee_i-\ee_{i+1}}^n\right] \left(f^{(\ss)}\right)^\kk
=\\=
e^{f(x)}\sum_{\{\kk|\kk\in\NN^\infty,\kk\cdot\ss=n+1\}} a_{\kk}^{n+1} \left(f^{(\ss)}\right)^\kk
\end{multline*}
This completes the proof.
\end{proof}

We are still not finished, though: we need to find an explicit formula for the coefficients. All we proved so far about them is that they can be determined recursively, hence uniquely, as this recursive determination makes them unique. The notation we are about to use is compatible with the notation used in \cite{S}. 

\begin{thm}
Let $\kk\in\NN^\infty$ and $n\in\NN$ so that $\ss\cdot \kk=n$. Then, $\ds a_\kk^n=\frac{1}{\prod_{i=1}^n k_i!}\frac{n!}{\prod_{i=1}^n (s_i!)^{k_i}}=
\frac{n!}{\kk!(\ss!)^\kk}$. For values of $\kk,n$ that do not satisfy the equality above, we set $\ds a_\kk^n=0$. 
\end{thm} 

\begin{proof}
All we need to show is that this definition is consistent with the recursive definition of the coefficients in the previous theorem. 

To begin with, observe that, assuming all quantities appearing are positive, 
\[a_{\kk-\ee_1}^n=\frac{1}{(k_1-1)!\prod_{i=2}^n k_i!}\frac{n!}{s_1^{k_1-1}\prod_{i=2}^n (s_i!)^{k_i}}=\frac{k_1s_1}{n+1}a_\kk^{n+1}\]
and
\begin{multline*}
a_{\kk+\ee_j-\ee_{j+1}}^n=\frac{1}{(k_j+1)!(k_{j+1}-1)!\prod_{i=1,i\neq j,j+1}^n k_i!}\frac{n!}{(s_j!)^{k_j+1}(s_{j+1}!)^{k_{j+1}-1}\prod_{i=1,i\neq j,j+1}^n (s_i!)^{k_i}}=\\=\frac{k_{j+1}}{k_j+1}\frac{s_{j+1}!}{s_j!}a_\kk^{n+1}=\frac{k_{j+1}}{k_j+1}s_{j+1}a_\kk^{n+1}
\end{multline*}

Now, observe that the equalities we have obtained are true even if not all the intermediate quantities are positive, as then the coefficients on both sides are 0, and so that the equalities become $0=0$. Then,
\[a_{\kk-\ee_1}^n+\sum_{i=1}^n (k_i+1)a_{\kk+\ee_i-\ee_{i+1}}^n=\frac{a_\kk^{n+1}}{n+1}\left(s_1k_1+\sum_{i=1}^n s_{i+1}k_{i+1}\right)=
a_\kk^{n+1}\frac{\ss\cdot\kk}{n+1}=a_\kk^{n+1}\]
Finally, we see immediately that $\ds a_1^1=\frac{1!}{1!1!}=1$, as it should be. This completes the proof. 
\end{proof}

\section{Number of summands}

How many are the solutions of the equation $\kk\cdot\ss=n\in\NN$ in $\NN^\infty$? They are as many as the different ways in which we can write $n$ as the sum of non-negative integers not greater than $n$, i.e. the number of partitions of $n$, as it is known in the literature (see \cite{A} for details). Although there exists an exact formula for this (see \cite{CG}), it is not practical at all; fortunately, simple bounds exist for this number. 

\section{Non-zero elements of the solutions}

How many non-zero elements does $\kk$ have for a given $n$? Obviously, at least one. On the other hand, the first time $s$ non-zero elements will appear in $\kk$ will be when they are all 1 and they occupy the lowest possible positions within the sequence, i.e. the first $s$ ones. Hence, $\kk\cdot\ss=n$ becomes $\ds n=1+\ldots+s=\frac{s(s+1)}{2}$. 

So, for a given $n$, what is the maximum number of non-zero elements there can be? As any nonzero element is as least 1, we need to solve $\ds \frac{s(s+1)}{2}\leq n$, or $s^2+s-2n\leq 0$, which leads to:
\[s\leq \left\lfloor\frac{\sqrt{1+8n}-1}{2}\right\rfloor\]  

\section{Conclusion}
Given the notation we stated in the beginning of the paper, we proved the identity:
\begin{equation}\left(e^{f(x)}\right)^{(n)}=e^{f(x)}\sum_{\{\kk|\kk\in\NN^\infty,\kk\cdot\ss=n\}} \frac{n!}{\kk!(\ss!)^\kk} \left(f^{(\ss)}\right)^\kk\end{equation}


\begin{thebibliography}{10}

\bibitem[1]{S}  R. P. Stanley. \textit{Enumerative Combinatorics, Vol. 1}\ \ Cambridge 1997

\bibitem[2]{A}  G. Andrews. \textit{The Theory of Partitions}\ \ Cambridge 1984

\bibitem[3]{CG} J. Conway, R. Guy. \textit{The Book of Numbers}\ \ Springer-Verlag 1996

\end{thebibliography}
\end{document}